\documentclass[12pt,reqno,fleqn]{amsart}  
\usepackage{tikz}
\usepackage{amsmath,amstext,amsthm,amssymb,amsxtra}
\usepackage[T1]{fontenc}

\usepackage{euscript}
\usepackage{mathrsfs}

\usepackage{mathtools}
\mathtoolsset{showonlyrefs,showmanualtags}

\usepackage{hyperref} 
\hypersetup{
	colorlinks=true,       
	linkcolor=blue,          
	citecolor=magenta,        
	filecolor=magenta,      
	urlcolor=cyan           
}

\usepackage[msc-links]{amsrefs}

\theoremstyle{plain} 
\newtheorem{lemma}{Lemma}[section]
\newtheorem{proposition}{Proposition}[section] 
\newtheorem{theorem}{Theorem} [section]

\theoremstyle{remark}
\newtheorem{remark}{Remark}

\def\supp{{\rm supp\,}}

\def\dist{{\rm dist}}

\def\ZL{\ensuremath{\mathcal L}}
\def\zN{\ensuremath{\mathcal N}}

\def\ZB{\ensuremath{\mathscr B}}
\def\ZC{\ensuremath{\mathscr C}}
\def\zB{\ensuremath{\mathbb B}}
\def\ZA{\ensuremath{\mathscr A}}

\def\ZI{\ensuremath{\mathbf 1}}
\def\ZN{\ensuremath{\mathbb N}}
\def\ZE{\ensuremath{\mathbb E}}

\def\kappa{\ensuremath{\mathcal K}}
\def\ZR{\ensuremath{\mathbb R}}

\def\ZF{\ensuremath{\mathscr F}}

\def\ZD{\ensuremath{\mathscr D}}

\theoremstyle{plain}
\newtheorem{definition}{Definition}[section]

\numberwithin{equation}{section}

\def\md#1#2\emd{\ifx0#1
	\begin{equation*} #2 \end{equation*}\fi  
	\ifx1#1\begin{equation}#2\end{equation}\fi   
	\ifx2#1\begin{align*}#2\end{align*}\fi   
	\ifx3#1\begin{align}#2\end{align}\fi    
	\ifx4#1\begin{gather*}#2\end{gather*}\fi  
	\ifx5#1\begin{gather}#2\end{gather}\fi   
	\ifx6#1\begin{multline*}#2\end{multline*}\fi  
	\ifx7#1\begin{multline}#2\end{multline}\fi  
	\ifx8#1\begin{multline*}\begin{split}#2\end{split}\end{multline*}\fi
	\ifx9#1\begin{multline}\begin{split}#2\end{split}\end{multline}\fi
}

\newcommand {\pro }[1]{Proposition \ref{#1}}

\newcommand {\df }[1]{Definition \ref{#1}}

\title[] {An equivalence principle for martingales }
\author[G. Karagulyan]{Grigori A. Karagulyan*}
\address{Institute of Mathematics of the National Academy of Sciences of the Republic of Armenia, Marshal Baghramian ave., 24/5, Yerevan, 0019, Armenia}
\email{g.karagulyan@gmail.com}

\author[L.  Khachatryan]{Linda A.  Khachatryan}
\address{Institute of Mathematics of the National Academy of Sciences of the Republic of Armenia, Marshal Baghramian ave., 24/5, Yerevan, 0019, Armenia} 
\email{linda@instmath.sci.am}

\author[Kh. Navoyan]{Khazhakanush V. Navoyan}
\address{Institute of Mathematics of the National Academy of Sciences of the Republic of Armenia, Marshal Baghramian ave., 24/5, Yerevan, 0019, Armenia} 
\email{khvnavoyan@gmail.com}

\thanks{Research was supported by the Higher Education and Science Committee of RA, in the frames of the
	research project 21AG-1A045}

\subjclass[2010]{60G42, 42C05, 42C10, 28D05}
\keywords{martingale difference, Haar system, non-overlapping polynomials, Weyl multiplier, Menshov-Rademacher theorem}
\begin{document}
	\begin{abstract}
	We show that any discrete-time filtration can be embedded into a dyadic filtration, realized via a non-classical measure-preserving map. This embedding makes it possible to reduce certain problems involving vector-valued martingales to the case of Haar martingales. Some applications of this result will also be discussed.
	\end{abstract}

	\maketitle  
	\section{Introduction}
	
	A sequence of $\sigma$-algebras $\{\ZF_n\}_{n\ge 1}$ in a probability space $(X,\ZA,\mu)$ is said to be a filtration if $\ZF_n\subset \ZF_{n+1}$. 
	A partition of a probability space $(X,\ZA,\mu)$ is a finite or infinite family of pairwise disjoint measurable sets $\{E_k\}\subset \ZA$ such that $\mu(E_k)>0$, $\cup_kE_k=X$. For partitions $\ZA$ and $\ZB$ we write $\ZA\prec \ZB$, if each element of $\ZA$ is a union of some elements of $\ZB$. A sequence of partitions $\{\ZF_n\}_{n\ge 1}$ is said to be a partitioning if
	\begin{equation}
		\ZF_n\prec \ZF_{n+1}, \quad n\ge 1, \text{ and } \lim_{n\to\infty}\max_{a\in \ZF_n}\mu(a)=0. 
	\end{equation}
	In particular, $\{\ZF_n\}_{n\ge 1}$ can be a partitioning on the probability space $[0,1)$ equipped with the Lebesgue measure. We refer to these as an interval partitioning if each $\ZF_n$ consists of intervals of the form $[a,b)$. If all those intervals are dyadic, then we call such a partitioning a dyadic partitioning. Given partitioning $\{\ZF_n\}_{n\ge 1}$ generates a filtration $\{\sigma(\ZF_n)\}_{n\ge 1}$, where $\sigma(\ZF_n)$ denotes the $\sigma$-algebra generated by the elements of the partition $\ZF_n$. A key result of the paper is an embedding theorem that allows a general filtration to be reflected into a dyadic one. Such a embedding is realized by a non-standard measure-preserving map between two probability spaces. Such a measure-preserving map is not a pointwise correspondence; rather, it acts on the measurable sets of the probability spaces. Namely,
	\begin{definition}
		Let $(X,\ZA,\mu)$ and $(Y,\ZB,\nu)$ be probability spaces. We say a function $\Theta:\ZA \to \ZB$ is measure-preserving set map (MPS-map), if
		\begin{align}
			&\text{1) }\mu(E)=\nu(\Theta(E)),\\
			&\text{2) }\mu(E\cup F)=\nu(\Theta(E)\cup \Theta(F)),\label{x43}
		\end{align}
		for any $E,F\in \ZA$. If the probability spaces $X$ and $Y$ coincide, then $\Theta$ is called a measure-preserving transformation (MPS-transformation).
	\end{definition} 
	We will examine in Section \ref{S2} the properties of such a measure-preserving map and explore its relationship to the classical notion of a measure-preserving transformation in the next section. In particular, we show that any MPS-map admits a natural extension to a map between the spaces of measurable functions on $X$ and $Y$ that preserves joint cumulative distribution functions (see \df{D1}). Specifically, for a simple function $f = \sum_j c_j \ZI_{a_j}$, its $\Theta$-image is defined by
	$\Theta(f)=\sum_jc_j\ZI_{\Theta(a_j)}$. For an arbitrary measurable function $f$ on $X$, the function $\Theta(f)$ is defined as the almost everywhere limit of $\Theta(s_n)$, where ${s_n}$ is a sequence of simple functions converging to $f$ almost everywhere (see Lemmas \ref{L1} and \ref{L6} for details).
	
	Introduce some preliminary notations. Define the distance of two partitions $\ZA$ and $\ZB$ as 
	\begin{equation}
		\dist(\ZA,\ZB)=\sum_{a\in \ZA\setminus \ZB}\mu(a)=\sum_{b\in \ZB\setminus \ZA}\mu(b).
	\end{equation}
	Note that if $\delta$ is the distance between two partitions, then they share common elements whose union has measure $1-\delta$. In particular, if the distance is zero, the two partitions coincide. For a partition $\ZF$ we denote by $\sigma(\ZF)$ the $\sigma$-algebra generated by the elements of $\ZF$. Then the conditional expectation of a function $f\in L^1(X)$ with respect to the $\sigma$-algebra $\sigma(\ZF)$ will be denoted by 
	\begin{equation}
		\ZE^\ZA f(x)=\sum_{a\in \ZA}\frac{1}{\mu(a)}\int_af\cdot \ZI_a(x),
	\end{equation} 
	where $\ZI_a$ is the indicator function of the set $a$. Given partitioning $\{\ZF_n\}_{n\ge 1}$ denote by $\ZF_\infty$ the $\sigma$-algebra generated by all the elements of partitions $\ZF_n$, i.e. $\ZF_\infty=\sigma\left(\cup_{n\ge 1} \ZF_n\right)$. 
	\begin{definition}
		Let $\{\ZA_n\}_{n\ge1}$ and $\{\ZB_n\}_{n\ge1}$ be partitionings on probability spaces $(X,\ZA,\mu)$ and $(Y,\ZB,\nu)$ respectively. The partitioning $\{\ZA_n\}_{n\ge1}$  is said to be embedded in $\{\ZB_n\}_{n\ge1}$ if there exists a MPS-map $\Theta:\ZA_\infty\to \ZB_\infty$ such that 
		\begin{equation}\label{x32}
			\Theta(\ZE^{\ZA_n}(\ZI_E))= \ZE^{\ZB_n}(\ZI_{\Theta(E)}) \text{ for all }E\in \cup_{j\ge 1}\ZA_j\text{ and }n\ge 1.
		\end{equation}
		The embedding of the partitioning $\{\ZA_n\}_{n\ge1}$ into $\{\ZB_n\}_{n\ge1}$ will be denoted by the symbol $\{\ZA_n\}_{n\ge1}\hookrightarrow\{\ZB_n\}_{n\ge1}$.	
	\end{definition}
	The following is the main result of the paper.
	\begin{theorem}\label{T1}
		Let $\varepsilon_n>0$, $n\ge 1$ and suppose $\{\ZA_n\}_{n\ge 1}$ and $\{\ZB_n\}_{n\ge 1}$ be partitionings on a probability space $X$. Then there exist sequence of integers $k_n$ and a partitioning $\{ \ZD_n\}_{n\ge 1}$ such that 
		\begin{align}
			&\ZD_n\subset \cup_{k\ge 1}\ZB_k,\label{x61}\\
			&\dist(\ZD_n,\ZB_{k_n})<\varepsilon_n,\label{x62}\\
			&\{\ZA_n\}_{n\ge 1}\hookrightarrow \{\ZD_n\}_{n\ge 1}.\label{x7}
		\end{align}
	\end{theorem}
	
	In Section \ref{S4} we will apply Theorem \ref{T1} to provide a framework for interpreting an arbitrary martingale (or martingale difference sequence) in terms of the Haar system (Theorem \ref{C1}). In Section \ref{S5} this equivalence principle is used to establish a logarithmic $L^p$-norm estimate for the maximal function of $N$ different transformations of a martingale. Finally, in Section \ref{S6} we deduce from Theorem \ref{T1} its multiple version.

	\section{Measure-preserving transformations }\label{S2}
	Let $(X,\ZA,\mu)$ be a probability space. By the classical definition, a measure-preserving transformation (MP-transformation) is a function $\phi: X\to X$ such that $\phi^{-1}(E)\in \ZA$ for every $E\in \ZA$ and $\mu(\phi^{-1}(E))=\mu(E)$. Observe that any MP-transformation $\phi: X\to X$ generates a MPS-transformation $\Theta: \ZA\to \ZA$ defined by $\Theta(E)=\phi^{-1}(E)$. Nevertheless, not all MPS-transformations can be represented by a MP-transformation. An example of such an MPS-transformation can be seen in the proof  of Theorem \ref{T1}. 
	
	\begin{lemma}\label{L1}
		If  $\Theta$ is a MPS-map  from a probability space $(X,\ZA,\mu)$ to $(Y,\ZB,\nu)$, then for any measurable sets $E,F\in \ZA$, $E_k\in \ZA$, $k=1,2,\ldots$,
		\begin{align}
			&\mu(E\setminus F)=\nu(\Theta(E)\setminus \Theta(F)),\label{x6}\\
			&\mu\left(\bigcup_k E_k\right)=\nu\left(\bigcup_k \Theta(E_k)\right),\label{x42}\\
			&\mu\left(\bigcap_k E_k\right)=\nu\left(\bigcap_k \Theta(E_k)\right).\label{a0}
		\end{align}
	\end{lemma}
	\begin{proof}
		Using the definition of MPS-map, we obtain
		\begin{align*}
			\mu(F)+\mu(E\setminus F)&=\mu(F\cup E)=\nu(\Theta(F)\cup\Theta(E))\\
			&=\nu(\Theta(F))+\nu(\Theta(E)\setminus\Theta(F))\\
			&=\mu(F)+\nu(\Theta(E)\setminus\Theta(F)).
		\end{align*}
		Thus \eqref{x6} follows. Similarly,
		\begin{align*}
			\mu(E)+\mu(F)-\mu(E\cap F)&=\mu(E\cup F)=\nu(\Theta(E)\cup \Theta(F)))\\
			&=\nu(\Theta(E))+ \nu(\Theta(F))-\nu(\Theta(E)\cap \Theta(F))\\
			&=\mu(E)+\mu(F)-\nu(\Theta(E)\cap \Theta(F)).
		\end{align*}
		This implies \eqref{a0} for two sets and so for a finite collection of measurable sets. Then passing to a limit we will get \eqref{a0} for any countable collection of sets. Similarly \eqref{x43} implies \eqref{x42}. 
	\end{proof}
	Throughout the paper, $\zB$ denotes a Banach space. Let $(X,\ZA,\mu)$ be a probability space. We denote by $L^0_{\zB}(X)$ the space of $\zB$-valued measurable functions on $X$. For an arbitrary family of measurable set $\ZF \subset \ZA$, we denote by $L^0_{\zB}(X, \ZF)$ the space of $\sigma(\ZF)$-measurable $\zB$-valued functions on $X$. Similarly, the notation $L^p_{\zB}(X, \ZF)$, $1\le p\le \infty$, denotes the Bochner-Lebesgue spaces of $\sigma(\ZF)$-measurable $\zB$-valued functions on $X$.
	
	Let $\Theta$ be a MPS-map from a probability space $(X,\ZA,\mu)$ to $(Y,\ZB,\nu)$. Then $\Theta$ induces a map that takes any indicator function $\ZI_E$, $E\in \ZA$, to $\ZI_{\Theta(E)}$. The following lemma gives a suitable extension of $\Theta$ to entire space $L^0_\zB(X)$ of $\zB$-valued measurable functions. 
	\begin{lemma}\label{L6}
		If  $\Theta$ is a MPS-map  from a probability space $(X,\ZA,\mu)$ to a probability space $(Y,\ZB,\nu)$, then it can be uniquely extended as an operator from $L^0_\zB(X)$ into $L^0_\zB(Y)$ such that  every sequence $\{f_k\}_{k\ge1}\subset L^0_\zB(X)$ is probabilistically equivalent to its image sequence $\{\Theta(f_k)\}_{k\ge 1}$.
	\end{lemma}
	\begin{proof}
		Let $s(x)$ be a simple function, i.e.
		\begin{equation}\label{a1}
			s(x)=\sum_{i=1}^m\alpha_i\ZI_{E_i}(x),
		\end{equation}
		where $E_i\in \ZA$, $i=1,2,\ldots,m $, are pairwise disjoint measurable sets and $\alpha_i\in \zB$. Then we define 
		\begin{equation}
			\Theta s(x)=\sum_{i=1}^m\alpha_i\ZI_{\Theta(E_i)}(x).
		\end{equation}
		For an arbitrary  $f\in L^0_\zB(X)$, one can find a sequence of simple functions $s_n(x)$, such that 
		\begin{equation}\label{a2}
			f(x)=\lim_{n\to\infty}s_n(x)\text{ a.e..}
		\end{equation}
		Using Lemma \ref{L1} and a standard argument one can observe that the sequences $s_n(x)$ and $\Theta s_n(x)$ are probabilistically equivalent. Consequently, $\Theta s_n(x)$ converges almost everywhere and the limit doesn't depend on the particular choice of simple functions $s_n$ in \eqref{a2}. We therefore define this limit as $\Theta f(x)$. It is also clear that any sequence $\{f_k\}_{k\ge1}\subset L^0_\zB(X)$ is probabilistically equivalent to its image sequence $\{\Theta(f_k)\}_{k\ge1}$.
	\end{proof}
	\begin{lemma}\label{L7}
		Let $\{\ZF_n\}_{n\ge 1}$ be a filtration on a probability space $(X,\ZF_\infty,\mu)$ and suppose that $\Theta: \cup_n\ZF_n\to  \ZF_\infty$ is a map, which restriction to each $\sigma$-algebra $\ZF_n$ is a MPS-map from $\ZF_n$ into $\ZF_\infty$. Then $\Theta$ can be extended as a MPS-transformation on $\ZF_\infty$.
	\end{lemma}
	\begin{proof}
		It is well known in measure theory that any set $E\in \ZF_\infty$ can be well-approximated by sets from $\cup_n\ZF_n$. In particular, for any set $E\in \ZF_\infty$ there exists a sequence $F_n\in \cup_n\ZF_n$ such that $\ZI_{F_n}\to \ZI_E$ a.e.. Since the sequence $\Theta(F_n)$ is probabilistically equivalent to the sequence $F_n$, we can say that $\ZI_{\Theta(F_n)}$ converges a.e. to an indicator function $\ZI_F$, $F\in \ZF_\infty$, and $F$ doesn't depend on the choice of the sequence $F_n$. Then we define $\Theta(E)=F$. In view of Lemma \ref{L1}, one can check that with this definition $\Theta$ becomes a MPS-transformation on $\ZF_\infty$.
	\end{proof}
	
	\section{Proof of Theorem \ref{T1}}
	\begin{lemma}\label{L9}
		Let $\{\ZA_n\}_{n\ge1}$ and $\{\ZB_n\}_{n\ge1}$ be partitionings on probability spaces $(X,\ZA,\mu)$ and $(Y,\ZB,\nu)$ respectively such that $\{\ZA_n\}_{n\ge1}\hookrightarrow\{\ZB_n\}_{n\ge1}$ and suppose that this embedding is realized by the MPS-transformation $\Theta$. Then,
		\begin{align}
			&1) f\in L^0_\zB(X,\ZA_n)\Rightarrow \Theta(f)\in L^0_\zB(X,\ZB_n),\label{x41}\\
			&2) \Theta(\ZE^{\ZA_n}(f))= \ZE^{\ZB_n}(\Theta(f)) \text{ for all }f\in L^1_\zB(X,\ZA_\infty)\text{ and }n\ge 1.\label{x31}
		\end{align}
	\end{lemma}
	\begin{proof}
		Observe that equation \eqref{x32} in the definition of embedding implies $\Theta(\sigma(\ZA_n))\subset \sigma(\ZB_n)$, which is the same as \eqref{x41} when $f=\ZI_E$ is an indicator function of a measurable set $E$. We also note that equation \eqref{x32} is a special case of \eqref{x31} when $f=\ZI_E$.
		The proofs of both relations \eqref{x41} and \eqref{x31} are based on a standard approximation argument using simple functions. 
	\end{proof}
	The following lemma follows directly from the definition of the embedding property of partitionings. 
	\begin{lemma}\label{L8}
		If partitionings $\{\ZA_n\}_{n\ge 1}$, $\{\ZB_n\}_{n\ge 1}$  and $\{\ZC_n\}_{n\ge 1}$ satisfy 
		\begin{equation}
			\{\ZA_n\}_{n\ge 1}\hookrightarrow \{\ZB_n\}_{n\ge 1},\quad \{\ZB_n\}_{n\ge 1}\hookrightarrow \{\ZC_n\}_{n\ge 1}
		\end{equation}
		realized by MPS-maps $\phi$ and $\psi$ respectively, then it holds the embedding
		\begin{equation*}
			\{\ZA_n\}_{n\ge 1}\hookrightarrow \{\ZC_n\}_{n\ge 1}
		\end{equation*}
		realized by $\phi\circ \psi$.
	\end{lemma}
	Let $a\subset [0,1)$ be an interval, and $a_k$ be a (finite or infinite) partition of $a$ into pairwise disjoint subintervals, 
	\begin{equation}\label{x16}
		a=\cup_ka_k,\quad a_k\cap a_j=\varnothing \text{ as } k\neq j.
	\end{equation}
	To this partition, we associate a Markov partition transformation $\phi:a\to a$, which maps each interval $a_k$ on $a$ linearly. Namely,
	\begin{equation}\label{x9}
		\phi(x)=\sum_k l_{a_k\to a}(x)\cdot \ZI_{a_k}(x),
	\end{equation}
	where $l_{a_k\to a}$ denotes the increasing linear function that maps $a_k$ onto $a$. Clearly, $\phi$ is a MP-transformation of the interval $a$.
	\begin{remark}\label{R1}
		Observe that from equation \eqref{x32} it follows that $\Theta(\ZA_n)\prec \ZB_n$. Thus to check equation \eqref{x32} we will need to show that
		\begin{equation}\label{x33}
			\frac{|\Theta(E)\cap \Theta(A)|}{|A|}=\frac{|\Theta(E)\cap B|}{|B|}
		\end{equation}
		for every choice of sets  
		\begin{equation}\label{x28}
			A\in \ZA_n,\quad B\in\ZB_n,\quad B\subset\Theta(A),\quad E\in \ZA_k.
		\end{equation}
	\end{remark}
	\begin{lemma}\label{L5}
		Let $\{\ZA_n\}_{n\ge 1}$ and $\{\ZC_n\}_{n\ge 1}$ be interval partitionings on $[0,1)$.  Then  for any $\varepsilon>0$ and for any fixed integer $m\ge 1$ there exist an interval partitioning $\{ \ZB_n\}_{n\ge 1}$ such that 
		\begin{align}
			&\ZB_m\subset \cup_{n\ge 1}\ZC_n,\label{x56}\\
			&\dist(\ZB_m,\ZC_{k_m})\le \varepsilon\text{ for an integer }k_m.\label{x57}
		\end{align}
		Besides, we have embedding
		\begin{equation}\label{x27}
			\{\ZA_n\}_{n\ge 1}\hookrightarrow \{\ZB_n\}_{n\ge 1}
		\end{equation}
		realized by a MPS-transformation $\theta_m$, satisfying
		\begin{equation}\label{x25}
			\theta_m(A)=A\in \ZB_k \text{ for all }A\in \ZA_k,\, k<m.
		\end{equation}
	\end{lemma}
	\begin{remark}\label{R2}
		One can check that from \eqref{x25} it follows that $\theta_m(\ZA_k)=\ZA_k\prec \ZB_k$ for all $k<m$.
	\end{remark}
	\begin{proof}
		For every interval $a\in \ZA_m$ we associate an arbitrarily chosen partition $\{a_j\}$ consisting of the elements of $\cup_{n\ge 1}\ZC_n$. Moreover, for sufficiently large $k_m$, most of the elements $\{a_j\}$ can be chosen from $\ZC_{k_m}$. More precisely, 
		\begin{equation}\label{x58}
			\sum_{a_j\in \ZC_{k_m}}|a_j|>|a|(1-\varepsilon) 
		\end{equation}
		and the integer $k_m$ is common for all elements $a\in \ZA_m$. 
		Let $\phi_{m,a}$ be the MP-transformation of the interval $a$, corresponding to its partition $\{a_j\}$ (see \eqref{x9}). Denote by $\phi_m$, the MP-transformation of $[0,1)$, such that $\phi_m(x)=\phi_{m,a}(x)$, $x\in a$, for all $a\in \ZA_m$. Define the MPS-transformation $\theta_m=\phi_m^{-1}$. Clearly, if $\{a_j\}$ is the defined partition of an element $a\in \ZA_m$, then for any measurable set $E\subset a$ we have
		\begin{equation}\label{x30}
			|\theta_m(E)\cap a_j|=\frac{|a_j|}{|a|}\cdot |E|.
		\end{equation}
		Set 
		\begin{equation}
			\ZB_n=\ZA_n\text{ if }n<m,
		\end{equation}
		while for $n\ge m$, $\ZB_n$ is defined as follows. For each $b\in \ZA_n$ there exists a unique $a\in \ZA_m$ such that  $b\subseteq a$. Moreover, $\theta_m(b)$ can be written as a union of the component-intervals $b_j=a_j\cap \theta_m(b)$ (see \eqref{x16}). We define the elements of $\ZB_n$ to be all such component-intervals $b_j$, where $b\in \ZA_n$. Hence the $\sigma$-algebras $\ZB_n$ for all $n\ge 1$ are defined and one can check that those form an interval partitioning. If $n=m$, then in the above notations we will have $b=a$ and $b_j=a_j$. Since $\{a_j\}\subset \cup_{n\ge 1}\ZC_n$, we get \eqref{x56}. Moreover, condition \eqref{x58} implies \eqref{x57}.  Observe that for any $a\in \ZA_m$ we have $\theta_m(a)=a\in \ZA_m$. Since each $b\in \ZA_n$ with $1\le n<m$ is a union of some elements of $\ZA_m$, we conclude that $\theta_m(E)=E$ for any $E\in \ZA_n$, $n<m$. This implies \eqref{x25}.
		To check embedding \eqref{x27}, we need to show equation \eqref{x33} (with $\Theta=\theta_m$) for every choice of the sets \eqref{x28}.
		Consider the following cases of integers $k$ and $n$ in \eqref{x28}:
		
		Case 1, $n<m$. Then from \eqref{x25} (see also Remark \ref{R2}) it follows that $\theta_m(A)=A=B$ and \eqref{x33} holds trivially.  
		
		Case 2,  $n\ge m$, $k<m$. Again by \eqref{x25} we have $\theta_m(E)=E$. If $A\cap E=\varnothing $, then $\theta_m(A)\cap \theta_m(E)=\varnothing $. Since $B\subset\theta_m(A)$ (\eqref{x28}), thus we get $B\cap \theta_m(E)=\varnothing $ and \eqref{x33} is trivial. If $A\cap E\neq\varnothing $, then $E\supset A$ so $\theta_m(E)\supset \theta_m(A)\supset B$ and we get \eqref{x33}.
		
		Case 3, $n\ge m$, $k\ge m$: There exists a unique $a\in \ZA_m$ such that $A\subset a$. Since $\theta_m(a)=a$ (see the definition of $\theta_m$), we get $B\subset \theta_m(A)\subset a$. If $E\cap a=\varnothing$ so $\theta_m(E)\cap a=\varnothing$, then equation \eqref{x33} is trivial. Hence we can assume that $E\subset a$. Let $\{a_j\}$ be the dyadic partition of $a$ considered in the definition of $\theta_m$. By the definition of $\ZB_n$ we have $B=a_{j_0}\cap\theta_m(A)$ for some $j_0$. Therefore,
		\begin{equation}
			B\cap \theta_m(E)=a_{j_0}\cap \theta_m(A)\cap \theta_m(E)=a_{j_0}\cap \theta_m(A\cap E).
		\end{equation}
		On the other hand, by \eqref{x30} we have
		\begin{align*}
			&|B|=|a_{j_0}\cap\theta_m(A)|=\frac{|a_{j_0}|}{|a|}\cdot |A|,\\ 
			&|B\cap \theta_m(E)|=|a_{j_0}\cap \theta_m(A\cap E)|=\frac{|a_{j_0}|}{|a|}\cdot |A\cap E|.
		\end{align*}
		Thus we obtain
		\begin{equation}
			\frac{	|\theta_m(E)\cap B|}{|B|}=\frac{|a_{j_0}|}{|a|\cdot |B|}\cdot |A\cap E|=\frac{|A\cap E|}{|A|}=\frac{|\theta_m(A)\cap \theta_m(E)|}{|A|}
		\end{equation}
		and so \eqref{x33}. This completes the proof of lemma.
	\end{proof}
	\begin{proof}[Proof of Theorem \ref{T1}]
		Without loss of generality we can suppose that our partitioning $\{\ZA_n\}_{n\ge 1}$ is an interval partitioning on $[0,1)$. Applying Lemma \ref{L5} consecutively for $m=1,2,\ldots$, we will get  interval partitionings $\{\ZA_n^{(m)}\}_{n\ge 1}$, $m=1,2,\ldots $ on $[0,1)$, satisfying
		\begin{align}
			&\ZA_m^{(m)}\subset \cup_{n\ge 1}\ZB_n,\quad m\ge 1,\label{x3}\\
			&\dist(\ZA_m^{(m)},\ZB_{k_m})< \varepsilon_m,\quad m\ge 1,\label{x59}\\
			&\{\ZA_n^{(m-1)}\}_{n\ge 1}\hookrightarrow \{\ZA_n^{(m)}\}_{n\ge 1},\quad m\ge 1,\label{x14}\\
			&\theta_m(A)=A\in \ZA_n^{(m)}\text{ for all }A\in \ZA_n^{(m-1)},\,1\le n<m, \label{x15}
		\end{align}
		where $\theta_m$ is the MPS-transformation realizing the embedding \eqref{x14} and $\ZA_n^{(0)}=\ZA_n$.
		Define $\ZD_n=\ZA_n^{(n)}$, $n\ge 1$. From \eqref{x15} (see also Remark \ref{R2}) it follows that
		\begin{equation}\label{x21}
			\ZD_{m-1}=\ZA_{m-1}^{(m-1)}=\theta_m(\ZA_{m-1}^{(m-1)})=\ZA_{m-1}^{(m)}\prec \ZA_{m}^{(m)}=\ZD_m.
		\end{equation}
		Thus $\{\ZD_m\}_{m\ge 1}$ provides a partitioning, and we obtain \eqref{x3}$\Rightarrow$\eqref{x61}, \eqref{x59}$\Rightarrow$\eqref{x62}.  Consider the MPS-transformation
		\begin{equation}\label{x11}
			\Theta_n=(\theta_n\circ\cdots\circ\theta_2\circ\theta_1).
		\end{equation}
		Observe that for any set $E\in \sigma(\ZA_n)$ we have 
		\begin{equation}\label{x12}
			\ZD_n=\ZA_n^{(n)}\ni\Theta_n(E)=\Theta_{n+1}(E)=\Theta_{n+2}(E)=\ldots.
		\end{equation}
		Indeed, using \eqref{x14} and Remark \ref{R1}, we can write
		\begin{align*}
			&E\subset \ZA_n=\ZA_n^{(0)}\Rightarrow \\
			&\theta_1(E)\in \ZA_n^{(1)} \Rightarrow \\
			&(\theta_2\circ \theta_1)(E)\in\ZA_n^{(2)}\Rightarrow\\
			&\vdots\\
			&(\theta_n\circ\ldots\circ\theta_2\circ\theta_1)(E)\in\ZA_n^{(n)}\Rightarrow \label{x20}\\
			&\Theta_n(E)\in \ZA_n^{(n)}=\ZD_n.
		\end{align*}
		Then from \eqref{x15} it follows that 
		\begin{equation}
			\Theta_n(E)=\Theta_{n+1}(E)=\Theta_{n+2}(E)=\ldots.
		\end{equation}
		Thus we obtain \eqref{x12}. Then for any set $E\in \sigma(\ZA_n)$ we can define 
		\begin{equation}\label{x23}
			\Theta(E)=\lim_{k\to\infty}\Theta_k(E)=\Theta_n(E).
		\end{equation}
		Hence $\Theta$ is defined for any $E\in \cup_{n\ge 1}\sigma(\ZA_n)$ and it gives a MPS-transformation induced on each $\sigma(\ZA_n)$. Applying Lemma \ref{L7} we can extend $\Theta$ up to a MPS-transformation on $\ZA_\infty$. Let us prove that $\Theta$ realizes the embedding $\{\ZA_n\}_{n\ge 1}\hookrightarrow \{\ZD_n\}_{n\ge 1}$. 
		Using Lemma \ref{L8}, we can say that $\Theta_m$ realizes the embedding
		\begin{equation}
			\{\ZA_n\}_{n\ge 1}\hookrightarrow \{\ZA_n^{(m)}\}_{n\ge 1}.
		\end{equation}
		In view of Remark \ref{R1} it means that
		\begin{equation}\label{x35}
			\frac{|\Theta_m(E)\cap \Theta_m(A)|}{|A|}=\frac{|\Theta_m(E)\cap B|}{|B|}
		\end{equation}
		whenever
		\begin{equation}\label{x36}
			A\in \ZA_n,\quad B\in \ZA_n^{(m)},\quad B\subset\Theta_m(A),\quad E\in \ZA_k.
		\end{equation}
		If $m\ge n$, then from \eqref{x15} it follows that $\ZA_n^{(m)}=\ZA_n^{(n)}=\ZD_n$ and we can replace $\ZA_n^{(m)}$ by $\ZD_n$ in \eqref{x36}. Then letting $m\to\infty$ and taking into account \eqref{x23}, we can also replace $\Theta_m$ by $\Theta$ in \eqref{x35} and \eqref{x36}. In such a new setting relation \eqref{x35} gives $\{\ZA_n\}_{n\ge 1}\hookrightarrow \{\ZD_n\}_{n\ge 1}$, completing the proof of theorem. 
	\end{proof}
	
	\section{An equivalence principle with Haar martingales}\label{S4}
	Let $\zB$ be a Banach space.  A martingale adapted to a partitioning $\{\ZF_n\}_{n\ge 1}$ is a sequence of functions $f_n:X\to \zB$ such that $f_n$ is $\sigma(\ZF_n)$-measurable and for the conditional expectation we have $\ZE^{\ZF_{n-1}}(f_n)=f_{n-1},\quad n\ge 2$. Given martingale $\{f_n\}_{n\ge 1} $ defines a martingale difference sequence $d_n=f_n-f_{n-1}$ (here $f_0=0$) so that $f_n=\sum_{k=1}^nd_k$, $n\ge 1$. Thus, the defining property of a martingale difference sequence is $\ZE^{\ZF_{n-1}}(d_n)=0,\,n\ge 2$.
	The simplest example of a martingale difference sequence is given by the Haar system ${h_n}$ in the form ${c_n h_n}$, where each coefficient $c_n \in \zB$. Recall the definition of Haar system on $[0,1)$. We have $h_1\equiv 1$. For an integer $n\ge 2$, with the representation $n=2^k+j$, $1\le j\le 2^k$, $k=0,1,\ldots$, denote 
	\begin{equation*}
		\Delta_n=\left[\frac{j-1}{2^k},\frac{j}{2^k}\right)
	\end{equation*}
	and let $\Delta_n^+$ and $\Delta_n^-$ be the left and the right halves of the interval $\Delta_n$. Then $n$'th Haar function for $n\ge 2$ is defined as 
	\begin{equation*}
		h_n=2^{k/2}\left(\ZI_{\Delta_n^+}-\ZI_{\Delta_n^-}\right).
	\end{equation*} 
	\begin{definition}
		A sequence of pairwise disjoint sets of positive integers $G_n\subset \ZN$, $n=1,2,\ldots$, is said to be locally increasing if for any integers $j$ and $k$ with
		\begin{equation}\label{x1}
			j\in G_n\text{ and }k\in G_{n+1}, 
		\end{equation}
		we have either $\Delta_k\subset \Delta_j$ or $\Delta_k\cap \Delta_j=\varnothing$. If $\max G_n<\min G_{n+1}$, then we say that the sequence $G_n$ is increasing.
	\end{definition}	
	Observe that if $G_n\subset \ZN$ is locally increasing sequence then for any coefficients $b_j\in \zB$ the functions 
	\begin{equation}\label{x29}
		f_n=\sum_{j\in G_n}b_jh_j
	\end{equation}
	form a martingale difference sequence. We also refer to such sequences of Haar sums as locally increasing.
	\begin{definition}\label{D1}
		Let $\zB$ be a Banach space and that $(X,\ZA,\mu)$ and $(Y,\ZB,\nu)$ are probability spaces. Sequences of measurable functions $f_k:X\to \zB$ and $g_k:Y\to \zB$, $k=1,2,\ldots$, are said to be probabilistically equivalent if they share the same joint cumulative distribution function, i.e. we have
		\begin{equation}\label{b8}
			\mu\left(\bigcap_{k=1}^nf_k^{-1}(E_k)\right)=\nu\left(\bigcap_{k=1}^ng_k^{-1}(E_k)\right)
		\end{equation}
		for every choice of Borel sets $E_k\subset \zB$, $k=1,2,\ldots,n$. Similarly we say that the sequences  of measurable sets $\{A_k\}\subset \ZA$ and $\{B_k\}\subset \ZB$ are probabilistically equivalent if the corresponding sequences of indicator functions $\{\ZI_{A_k}\}$ and $\{\ZI_{B_k}\}$ are probabilistically equivalent.
	\end{definition}
	\begin{theorem}\label{C1}
		Let $\zB$ be a separable Banach space and $\varepsilon_n\to 0$.  Then any $\zB$-valued atomic martingale difference sequence $\{d_n\}_{n\ge 1}$ is probabilistically equivalent to a sequence of $\zB$-valued Haar sums
		\eqref{x29}, where the sequence $G_n$ is locally increasing and $b_j\in \zB$. Moreover, there are subsets $G'_n\subset G_n$ such that $G'_n$ is increasing and 
		\begin{equation}\label{x40}
			\left|\left\{\sum_{j\in G_n}b_jh_j\neq \sum_{j\in G'_n}b_jh_j\right\}\right|<\varepsilon_n.
		\end{equation}
	\end{theorem}
	For $f\in L^1_\zB[0,1]$ we consider the corresponding vector-valued Haar series
	\begin{equation}\label{x22}
		\sum_{n=1}^\infty (f,h_n)h_n(x),\text{ where }(f,h_n)=\int_0^1f\cdot h_n.
	\end{equation} 
	The following standard properties of \eqref{x22} can be derived by applying the same argument as in the classical case, where $\zB=\ZR$.
	\begin{proposition}\label{P0}
		Let $\zB$ be an arbitrary Banach space. Then for any $f\in C_\zB[0,1]$ the series \eqref{x22} uniformly converges to $f$.
	\end{proposition}
	\begin{proposition}\label{P1}
		If $\zB$ is a separable Banach space, then for any $f\in L^p_\zB[0,1]$, $1\le p<\infty$, the series \eqref{x22} converges in $L^p$-norm and almost everywhere.
	\end{proposition}
	\begin{definition}
		Let $\ZF$ be partition of a probability space $(X,\ZF,\mu)$ and denote
		\begin{equation}
			m(\ZF)=\max_{a\in \ZF}\mu(a)
		\end{equation}
		We say that an element $a\in \ZF$ is maximal, if $\mu(a)=m(\ZA)$. A partition $\ZA$ is said to be $\varepsilon$-uniform, $\varepsilon>0$, if 
		\begin{equation*}
			\sum_{a\in \ZF),\, \mu(a)<m(\ZF)}\mu(a)<\varepsilon.
		\end{equation*}
	\end{definition}
	
	\begin{lemma}\label{L4}
		Let $\{\ZF_n\}_{n\ge 1}$ be a dyadic partitioning on $[0,1)$. Then any $\zB$-valued martingale difference $\{f_n\}_{n\ge 1}$ adapted to $\{\ZF_n\}_{n\ge 1}$ admits a Haar representation
		\begin{equation}\label{x39}
			f_n=\sum_{j\in G_n}(f_n,h_j)h_j,
		\end{equation}
		where the sets $G_n\subset \ZN$, $n=1,2,\ldots$, form a locally increasing sequence of positive integers. Moreover, if in addition $\ZF_n$ is $\varepsilon_n$-uniform, $\varepsilon_n>0$, then there exists an increasing sequence $G'_n$ such that $G'_n\subset G_n$ and
		\begin{equation}\label{x37}
			\left|\supp\left(\sum_{j\in G_n\setminus G'_n}(f_n,h_j)h_j\right)\right|<\varepsilon_n.
		\end{equation}
	\end{lemma}
	\begin{proof}
		Denote 
		\begin{align*}
			&G_1=\{k\in \ZN:\, \Delta_k\not\subseteq I \text{ for any } I\in \ZF_1\},\\
			&G_n=\{k\in \ZN:\, \Delta_k\not\subseteq I \text{ for any } I\in \ZF_n \text{ and }\\
			&\qquad\qquad \qquad \quad \Delta_k\subseteq I\text{ for some } I\in \ZF_{n-1}\},\quad n\ge 2.
		\end{align*}
		One can check that the sequence $\{G_n\}_{n\ge 1}$ is locally increasing. Then, in view of \pro{P1}, it remains to show that 
		\begin{equation}
			(f_n,h_k)=0\text{ for all }k\not\in G_n. 
		\end{equation}
		Indeed, if $k\not\in G_n$, then one of the following conditions hold:
		\begin{align}
			&1)\,\Delta_k\subseteq I \text{ for some } I\in \ZF_n,\\
			&2)\, \Delta_k\not\subseteq I\text{ for every } I\in \ZF_{n-1}.
		\end{align}
		Suppose we have 1), i.e. there is $I\in \ZF_n$ such that $\Delta_k\subseteq I$.  Since $f_n$ is constant on $I$, and therefore on $\Delta_k$, it follows that $(f_n,h_k)=0$. If a $\Delta_k$ satisfies 2), then we get that both $\Delta_k^+$ and $\Delta_k^-$ are unions of some elements of $\ZF_{n-1}$, since each interval $I\in \ZF_{n-1}$ itself is dyadic. On the other hand, by the martingale-difference property of $f_n$ we have $(f_n,\ZI_I)=0$ for any $I\in \ZF_{n-1}$. This implies $(f_n, \ZI_{\Delta_k^\pm})=0$, and therefore, $(f_n,h_k)=0$. If $\ZF_n$ is $\varepsilon_n$-uniform, then we have 
		\begin{equation}\label{x34}
			\sum_{a\in\ZF_n,\, \mu(a)<m(\ZF_n)}\mu(a)<\varepsilon_n.
		\end{equation}
		Denote
		\begin{align}\label{x38}
			G'_n=\{k\in G_n:\, &\text{there exists a maximal element } I\in \ZF_n\\
			&\text{ such that }\Delta_k \supsetneq I\}.
		\end{align}
		Observe that 
		\begin{equation}
			\Delta_k\subseteq \bigcup_{a\in\ZF_n,\, \mu(a)< m(\ZF_n)}a\text{ if }k\in G_n\setminus G'_n.
		\end{equation}
		Combining this with \eqref{x34}, we obtain \eqref{x37}. 
		Now consider integers $k\in G'_{n-1}$ and $j\in G'_{n}$. By the definitions of $G'_{n-1}$ (\eqref{x38}) we have $|\Delta_k|>m(\ZF_{n-1})$. Since $j\in G'_n\subset G_n$, by the definition of $G_n$ there exists an atom $I\in \ZF_{n-1}$ such that $\Delta_j\subseteq I$ and so $|\Delta_j|\le m(\ZF_{n-1})$. Thus we get $|\Delta_j|<|\Delta_k|$, which clearly implies $j>k$. Hence $\{G'_n\}$ is an increasing sequence of integer sets.
	\end{proof}
	\begin{proof}[Proof of Theorem \ref{C1}]
		Let $\{d_n\}_{n\ge 1}$ be a $\zB$-valued martingale difference sequence adapted to a partitioning $\{\ZF_n\}_{n\ge 1}$. Applying Theorem \ref{T1}, we find a dyadic  partitioning $\{\ZD_n\}_{n\ge 1}$ such that each $\ZD_n$ is $\varepsilon_n$-uniform and $\{\ZF_n\}_{n\ge 1}\hookrightarrow \{\ZD_n\}_{n\ge 1}$. Let $\Theta$ be the MPS-map realizing this embedding. Using the property $\Theta(\ZF_n)\subset \ZD_n$ (see Remark \ref{R1}) and applying Lemma \ref{L6}, we conclude that the sequence $\{f_n=\Theta(d_n)\}_{n\ge 1}$ is probabilistically equivalent to $\{d_n\}_{n\ge 1}$ and is adapted to the partitioning $\{\ZD_n\}_{n\ge 1}$. On the other hand according to Lemma \ref{L9} we have
		\begin{equation*}
			\ZE^{\ZD_{n-1}}(\Theta(d_n))=\Theta(\ZE^{\ZA_{n-1}}(d_n)).
		\end{equation*}
		Since $\{f_n\}_{n\ge 1}$ is a martingale difference, we have $\Theta(\ZE^{\ZA_{n-1}}(d_n))=0$. Therefore, 
		\begin{equation*}
			\ZE^{\ZD_{n-1}}(f_n)=\ZE^{\ZD_{n-1}}(\Theta(d_n))=0, 
		\end{equation*}
		i.e. $\{f_n\}_{n\ge 1}$ is also a martingale difference adapted to the dyadic partitioning $\{\ZD_n\}_{n\ge 1}$. Then applying Lemma \ref{L4}, we get the representation \eqref{x39},
		where $\{G_n\}$ is a locally increasing sequence of integer sets. Using also \eqref{x37}, we get \eqref{x40}.
		This completes the proof of theorem.
	\end{proof}
	
	\section{Martingale transforms}\label{S5}
	Vector-valued Haar martingales represent the simplest models of general vector-valued martingales. Owing to Theorem \ref{C1}, the verification of various martingale inequalities in the vector-valued setting can be reduced to the dyadic Haar case, from which the corresponding results for arbitrary vector-valued martingales follow. Important examples of such inequalities are those arising in the study of martingale transforms. 
	It is of interest Banach spaces $\zB$ for which classical theorems on real-valued martingales can be extended for $\zB$-valued martingales. A Banach space $\zB$ is said to have a UMD-property if there exists a positive constant $c_p$ such that the inequality 
	\begin{equation}\label{x51}
		\left\|\varepsilon_1d_1+\cdots+\varepsilon_nd_n\right\|_{p}\le c_p \left\|d_1+\cdots+d_n\right\|_{p},\quad 1<p<\infty,
	\end{equation}
	holds for all $\zB$-valued martingale difference sequences $\{d_j\}_{j\ge 1}$, for all so called predictable sequences $\varepsilon_j=\pm 1$ and all $n\ge 1$.  A geometric characterization of  UMD-spaces was given by Burkholder \cite{Burk1}. It is well known that $\ZR$ is a UMD-space. For Haar scalar-martingales this was proved by Paley \cite{Paley}. The inequality for the general scalar-martingales was established by  Burkholder \cite{Burk2}. Moreover, it was found in \cite{Burk2} the best constant $c_p$ that can be in \eqref{x51}, where $\{d_k\}$ are scalar-valued martingales. In view of Theorem \ref{C1}, problems of this type may, without loss of generality, be restricted to Haar martingales, taking $d_k=c_kh_k$ with $c_k\in \zB$. In particular, if one seeks the best constant in \eqref{x51} for a given Banach space $\zB$, then the same constant is obtained by considering only $\zB$-valued Haar martingale difference sequences instead of general ones. Indeed, by approximation, the optimal constant 
	$c_p$ remains unchanged when passing from arbitrary martingale differences to atomic ones. Then, applying Theorem \ref{C1}, for atomic martingale differences 
	in \eqref{x51} we find a sequence
	\begin{equation}\label{x53}
		f_k=\sum_{j\in G_k}b_jh_j,\quad k=1,2,\ldots,n,
	\end{equation}
	which is probabilistically equivalent to the $d_1,\ldots,d_n$. Then the inequality \eqref{x51} can be equivalently written as 
	\begin{equation}\label{x52}
		\left\|\varepsilon_1f_1+\cdots+\varepsilon_nf_n\right\|_{p}\le c_p \left\|f_1+\cdots+f_n\right\|_{p}.
	\end{equation}
	On the other hand, both $f_1+\cdots+f_n$ and $\varepsilon_1f_1+\cdots+\varepsilon_nf_n$ become Haar sums after substitutions \eqref{x53}, where the latter is simply the martingale transformation of the former. Thus we obtain the following result:
	\begin{theorem}\label{T2}
		Let $\zB$ be a UMD-space. If \eqref{x51} holds for all $\zB$-valued Haar martingales and all choices $\varepsilon_k=\pm 1$, then it holds for all $\zB$-valued martingale difference sequences, with the same constant $c_p$.
	\end{theorem}
	The inequalities in \eqref{x51}  may also be considered with $\varepsilon_k$ taking values in $[0,1]$ or in $\{0,1\}$. The statement of Theorem \ref{T2} remains valid for such choices of predictable $\varepsilon_k$.  In particular,  the case of $\varepsilon_k\in \{0,1\}$ was studied by Choi \cite{Choi}, who determined the optimal constant in \eqref{x51} for scalar-valued martingale difference sequences. 
	
	Now let $\{\ZF_n\}_{n\ge 1}$ be an arbitrary filtration in a probability space $(X,\ZF_\infty,\mu)$ and suppose that $f=\{f_n\}_{n\ge 1}$ is a scalar-valued martingale adapted to $\{\ZF_n\}_{n\ge 1}$ with assumption that $f_0=0$. For $1\le p\le \infty$ denote 
	\begin{equation*}
		\|f\|_p=\sup_{n\ge 1} \|f_n\|_p,\quad f^*=\sup_{n\ge 1}|f_n|.
	\end{equation*}
	The martingale $L^p$-space consists of all martingales with finite norm $\|f\|_p<\infty$. It is well known that for $1<p<\infty$, the martingale $L^p$-space is isomorphic to $L^p(X,\ZF_\infty)$.  Moreover, if $f=\{f_n\}_{n\ge 1}$ belongs to a martingale $L^p$-space, then $f_n$ converges both almost everywhere and in the $L^p$-norm. In this case, by convention, when writing a martingale $f=\{f_n\}_{n\ge 1}$ we identify $f$ with the limit function of $f_n$.  If $\Theta=\{\theta_n\}_{n\ge 1}$ is a sequence of real valued functions with $|\theta_n|\le 1$, then we say $\Theta$ is predictable if $\theta_n$ is $\ZF_{n-1}$-measurable.  The sequence $g=\{g_n\}_{n\ge 1}$ defined by
	\begin{equation}
		g_n=\sum_{k= 1}^n \theta_k d_k,\quad d_k=f_k-f_{k-1},
	\end{equation}
	is also a martingale and it is called the martingale transform of $f$ by $\Theta$. 
	Let denote by $T_\Theta$ the operator assigning any martingale $f$ to its $\Theta$-transformation $g$.
	From the inequality \eqref{x51} it follows that $T_\Theta$ is a bounded operator on $L^p$, $1<p<\infty$.  Moreover, we have the inequality
	\begin{equation}
		\|T_\Theta^*(f)\|_p\le c_p\|f\|_p,
	\end{equation}
	which is the combination of Burkholder's \cite{Burk1} and Doob's \cite{Doob} inequalities. For a finite family of predictables $\Theta_k=\{\theta^{(k)}_n\}_{n\ge 1}$, $k=1,2,\ldots,N$, consider the maximal operator
	\begin{equation}\label{x54}
		T(f)=\max_{1\le k\le N}T_{\Theta_k}^*(f).
	\end{equation}
	Applying Theorem \ref{C1}, we establish the following inequality, extending the $p=2$ case proved in \cite{Kar1}.  Moreover, in \cite{Kar1}, martingale transforms are considered in the case where the predictable sequences consist of constant functions.
	\begin{theorem}
		Let $\{\ZF_n\}_{n\ge 1}$ be a filtration in a probability space $(X,\ZF_\infty,\mu)$. Then operator \eqref{x54} satisfies
		\begin{equation*}
			\|T(f)\|_p\le c_p\sqrt{\log (N+1) }\cdot \|f\|_p, \quad 1<p<\infty,
		\end{equation*}
		for any $f\in L^p(X,\ZF_\infty)$.
	\end{theorem}	
	\begin{proof}
		By Theorem \ref{C1} it is enough to prove the theorem for Haar martingales. Let $\sum_{k=1}^\infty a_kh_k$ be the Fourier-Haar series of a function $f\in L^p(0,1)$. Recall the definitions of the maximal function and the square function 
		\begin{equation*}
			f^*(x)=\sup_{n\ge 1}\left|\sum_{k=1}^na_kh_k(x)\right|,\quad Sf(x)=\left(\sum_{k=1}^\infty a_k^2h_k^2(x)\right)^{1/2}.
		\end{equation*}
		It is well known the boundedness of both operators on $L^p$, $1<p<\infty$  (see for example \cite{KaSa}, ch. 3, Theorems 4 and 9). A key point in the proof of Theorem \ref{T1} is the following good-$\lambda$ inequality due to Chang-Wilson-Wolff  (see \cite{CWW}, Corollary 3.1):
		\begin{equation}\label{CWW}
			|\{x\in[0,1):\,f^*(x)>\lambda,\, Sf(x)<\varepsilon\lambda\}|\\\lesssim\exp\left(-\frac{c}{\varepsilon^2}\right)|\{f^*(x)>\lambda/2\}|,\,\lambda>0,\,0<\varepsilon <1.
		\end{equation}
		Denote
		\begin{equation*}
			f_k=T_{\Theta_k}(f)=\sum_{j=1}^\infty \theta_j^{(k)}a_j h_j,\quad k=1,2,\ldots,N.
		\end{equation*}
		We have $Sf_k(x)\le Sf(x)$, $k=1,2,\ldots,N$. Thus, applying inequality \eqref{CWW} with $\varepsilon_N=(c/\ln N)^{1/2}$, we obtain 
		\begin{align}\label{b2}
			|\{f_k^*>\lambda,\, S(f)\le\varepsilon_N\lambda\}|\lesssim\exp\left(-\frac{c}{\varepsilon_N^2}\right)|\{f_k^*>\lambda/2\}|.
		\end{align}
		For 
		\begin{equation*}
			T(f)=\max_{1\le k\le N}T_{\Theta_k}^*(f)=\max_{1\le k\le N}f_k^*
		\end{equation*}
		we obviously have
		\begin{align}
			\{ T(f)>\lambda\}&\subset \{T(f)>\lambda,\, S(f)\le \varepsilon_N\lambda\}\\
			&\cup  \{S(f)> \varepsilon_N\lambda\}=A(\lambda)\cup B(\lambda),
		\end{align}
		and thus
		\begin{equation}
			\|T(f)\|_p^p\le p\int_0^\infty\lambda^{p-1} |A(\lambda)|d\lambda+p\int_0^\infty\lambda^{p-1} |B(\lambda)|d\lambda.
		\end{equation}
		From \eqref{b2} it follows that
		\begin{align}
			p\int_0^\infty\lambda^{p-1}|A(\lambda)|d\lambda&\le p\sum_{k=1}^N\int_0^\infty\lambda^{p-1} |\{ f_k^*>\lambda,\, S(f)\le \varepsilon_N\lambda\}|d\lambda\\
			&\le \exp\left(-\frac{c}{\varepsilon_N^2}\right)\sum_{k=1}^Np\int_0^\infty\lambda^{p-1} |\{ f_k^*>\lambda/2\}|d\lambda\\
			&\lesssim\frac{1}{N}\sum_{k=1}^N\|f_k^*\|_p^p\\
			&\lesssim \frac{1}{N}\sum_{k=1}^N\| f_k\|_p^p\\
			&\lesssim \|f\|_p^p.
		\end{align}
		Combining this and 
		\begin{align*}
			p\int_0^\infty\lambda^{p-1}|B(\lambda)|d\lambda&=\varepsilon_N^{-p}\|S(f)\|_p^p\lesssim (\log N )^{p/2}\cdot \|f\|_p^p,
		\end{align*}
		we get
		\begin{equation}\label{b4}
			\|T(f)\|_p\lesssim \sqrt{\log N}\cdot \|f\|_p
		\end{equation} 
		which implies the theorem for Haar martingales, and therefore for arbitrary martingales as well.
	\end{proof}
	\section{Remark on multiple martingales}\label{S6}
	
	We consider the probability space $(Q_d,\ZL_d,m_d)$ of Lebesgue measure on $d$-dimensional unit cube $Q_d=[0,1)^d$.  Let $\{\ZA_n^{(k)},\,n\in \ZN \}$, $k=1,2,\ldots ,d$, be partitionings of $[0,1)$. The Cartesian product 
	\begin{equation}\label{I1}
		\ZA_{\bf n}=\ZA_{n_1,\ldots,n_d}=\otimes_{k=1}^d\ZA_{n_k}^{(k)},\quad {\bf n}=(n_1,\ldots,n_d)\in \ZN^d.
	\end{equation}
	determines a multiple partitioning of the unite cube $Q_d$. For a given sequence of integers $\zN=\{\nu_i\in \ZN:\, \nu_1<\nu_2<\ldots \}$ consider the multiple sequence of partitioning
	\begin{equation}\label{i5}
		\tilde \ZA_{\bf n}=\ZA_{\nu_{n_1},\ldots, \nu_{n_d}},\quad {\bf n}=(n_1,\ldots,n_d)\in\ZN^d.
	\end{equation}
	Partitioning  \eqref{i5} is said to be a sub-partitioning of \eqref{I1}. The following theorem is an improved version of the analogous result from \cite{Kar2}. It is obtained by a direct application of Theorem \ref{T1}.
	\begin{theorem}\label{T5}
		Let $\{\ZA_{\bf n},\, {\bf n}\in \ZN^d\}$ and $\{\ZB_{\bf n},\, {\bf n}\in \ZN^d\}$ be arbitrary multiple interval partitionings of the unit cube $Q_d$. Then for any sequence of numbers $\varepsilon_k>0$, $k=1,2,\ldots$, there exist a multiple partitioning $\ZD_{\bf n}$, a MPS-transformation $\Theta$, and a sub-partitioning $\{\tilde \ZB_{\bf n}\}$ of $\{\ZB_{\bf n}\}$, such that
		\begin{align}
			&\ZD_{\bf n}\subset \cup_{\bf k\ge 1}\ZB_{\bf k},\label{x81}\\
			&\dist(\ZD_{\bf n},\tilde \ZB_{\bf n})<\varepsilon_m,\quad \min_{1\le i\le d}n_i\ge m,\label{x82}\\
			&\left(\Theta\circ\ZE^{\ZA_{\bf n}}\right)f(x)=\left(\ZE^{\ZD_{\bf n}}\circ\Theta\right) f(x),\quad {\bf n}\in \ZN^d,\label{x83}
		\end{align}
		for all $f\in L^1(Q_d)$.
	\end{theorem}
	If $T$ is a bounded linear operator on $L^1[0,1)$, then we denote by $T_k$, $1\le k\le d$, the operator
	\begin{equation}
		T_k:L^1( Q_d)\to L^1(Q_d) 
	\end{equation}
	defined by
	\begin{equation}
		T_kf(x_1,\ldots,x_d)=Tf(x_1,\ldots,x_{k-1},\cdot ,x_{k+1},\ldots,x_d).
	\end{equation}
	The following lemmas are proved in \cite{Kar1}.
	\begin{lemma}[\cite{Kar1}]\label{L}
		If $T$ is a bounded linear operator on $L^1[0,1)$, then $T_k$ determines a bounded operator on $L^1(Q_d)$.
	\end{lemma}
	\begin{lemma}[\cite{Kar1}]\label{L0}
		If $T$ and $P$ are bounded linear operators on $L^1[0,1)$, then
		\begin{equation}
			P_n\circ T_m=T_m\circ P_n,\quad n\neq m,\quad 1\le n,m\le d.
		\end{equation}
	\end{lemma}
	\begin{proof}[Proof of Theorem \ref{T5}]
		Applying Theorem \ref{T1} for each $k=1,2,\ldots,d$ we find a partitioning  $\{\ZD_n^{(k)}\}_{n\ge 1}$ of the interval $[0,1)$ such that 
		\begin{align}
			&\ZD_n^{(k)}\subset \cup_{m\ge 1}\ZB_m^{(k)},\label{x71}\\
			&\dist(\ZD_m^{(k)},\ZB_{\nu_m}^{(k)})<\varepsilon_m/d,\label{x72}\\
			&\{\ZA_n^{(k)}\}_{n\ge 1}\hookrightarrow \{\ZD_n^{(k)}\}_{n\ge 1}.\label{x73}
		\end{align}
		Clearly, \eqref{x71} $\Rightarrow$ \eqref{x81}, and \eqref{x72} $\Rightarrow$ \eqref{x82}. Let $\theta_k$ be the MPS-transformation of $[0,1)$ realizing \eqref{x73}. By Lemma \ref{L6} it can be extended as a linear operator on $L^1[0,1)$. Then by Lemma \ref{L9} we can write
		\begin{equation*}
			\left(\theta_k\circ\ZE^{\ZA_{n_k}^{(k)}}\right)(f)=\left(\ZE^{\ZD_{n_k}^{(k)}}\circ\theta_k\right) (f),\text{ for all } f\in L^1[0,1).
		\end{equation*}
		Applying $\theta_k$ along $k$th coordinate, by Lemma \ref{L} we can consider $\theta_k$ as a linear bounded operator on $L^1(Q_d)$. Thus $\Theta=\theta_1\circ\theta_2\circ\cdots\circ\theta_d$ is also a linear bounded operator on $L^1(Q_d)$. We also have
		\begin{equation*}
			\ZE^{\ZD_{\bf n}}=\ZE^{\ZD_{n_1}^{(1)}}\circ\cdots\circ \ZE^{\ZD_{n_d}^{(d)}},\quad \ZE^{\ZA_{\bf n}}=\ZE^{\ZA_{n_1}^{(1)}}\circ\cdots\circ \ZE^{\ZA_{n_d}^{(d)}}
		\end{equation*}
		Applying Lemma \ref{L0}, we can write 
		\begin{align*}
			\Theta\circ\ZE^{\ZA_{\bf n}}&=\theta_1\circ\cdots\circ\theta_d\circ \ZE^{\ZA_{n_1}^{(1)}}\circ\cdots\circ \ZE^{\ZA_{n_d}^{(d)}}\\
			&=\left(\theta_1\circ \ZE^{\ZA_{n_1}^{(1)}}\right)\circ\cdots\circ\left(\theta_d\circ \ZE^{\ZA_{n_d}^{(d)}}\right)\\
			&=\left( \ZE^{\ZD_{n_1}^{(1)}}\circ\theta_1\right)\circ\cdots\circ\left(\ZE^{\ZD_{n_d}^{(d)}}\circ\theta_d \right)\\
			&=\ZE^{\ZD_{n_1}^{(1)}}\circ\cdots\circ \ZE^{\ZD_{n_d}^{(d)}}\circ\theta_1\circ\cdots\circ\theta_d\\
			&=\ZE^{\ZD_{\bf n}}\circ \Theta,
		\end{align*}
		which gives \eqref{x83}.
	\end{proof}
	
	\section*{Acknowledgments}
	Research was supported by the Higher Education and Science Committee of RA, in the frames of the research project 21AG-1A045
	
	\noindent
	\textbf{Conflicts of Interest}: The authors declare no conflicts of interest.
	
	\noindent
	\textbf{Data Availability Statement}: Data sharing is not applicable to this article as no data sets were generated or analyzed during the current study.
	
	\bibliographystyle{plain}
	

\end{document}